\documentclass[10pt]{article}
\usepackage{amsmath}
\usepackage{amssymb}
\usepackage{epsf}
\usepackage[all]{xy}
\usepackage[francais,english]{babel}

\newtheorem{theoreme}{Theorem}[section]
\newtheorem{propoe}[theoreme]{Proposition}
\newtheorem{coroe}[theoreme]{Corollary}
\newtheorem{dfne}[theoreme]{Definition}
\newtheorem{propdefe}[theoreme]{Definition-Proposition}
\newtheorem{rquee}[theoreme]{Remark}
\newtheorem{lemmae}[theoreme]{Lemma}
\newtheorem{csqcee}[theoreme]{Consequence}
\newtheorem{exple}[theoreme]{Example}
\newtheorem{proptee}{Property}
\newtheorem{notationse}[theoreme]{Notation}

\newenvironment{rmk}{\begin{rquee} \normalfont}{\end{rquee}}
\newenvironment{theo}{\begin{theoreme} \sf}{\end{theoreme}}
\newenvironment{prop}{\begin{propoe} \sf}{\end{propoe}}
\newenvironment{cor}{\begin{coroe} \sf}{\end{coroe}}

\newenvironment{lemma}{\begin{lemmae} \sf}{\end{lemmae}}

\newenvironment{proof}{{\flushleft{\bf Proof: }} \rm}{\Qed \\}
\newenvironment{prooflemma}{{\flushleft{\bf Proof of Lemma: }} \rm}{\qed \\}

\def\qed{\hfill \mbox{$\square$}}
\def\Qed{\hfill \mbox{$\blacksquare$}}

\def\ot{\otimes}
\def\gl{\Gamma_{\! \! \Lambda}}
\def\gln{\Gamma_{\! \! \Lambda _n}}

\def\lan{\Lambda _n}

\def\uq{$u_q^+(\frak{sl_2})$}
\newcommand{\topm}{\mathrm{top}}

\newcommand{\Ext}{\mathrm{Ext}}

\newcommand{\Tot}{\mathrm{Tot}}
\newcommand{\Hh}{\mathrm{HH}}
\newcommand{\Hc}{\mathrm{HC}}
\newcommand{\Hom}{\mathrm{Hom}}

\newcommand{\K}{\mathrm{K}}

\title{\textsc{Hochschild and cyclic homology of a family of Auslander algebras}}

\author{Rachel Taillefer\thanks{Supported by a Lavoisier grant from
    the French Ministry of Foreign Affairs.}}

\date{}

\begin{document}

\maketitle


\selectlanguage{english}
\begin{abstract} In this paper, we compute  the Hochschild and  cyclic
  homologies of the Auslander algebras of the Taft algebras. We also describe the first Chern character for the Taft
  algebras and for their Auslander algebras.
\end{abstract}

\paragraph{2000 Mathematics Subject Classification:} 16E20,  16E40,  19D55.
\paragraph{Keywords:} Hochschild homology, Cyclic homology, Auslander algebras, Chern characters.

\section{Introduction}

The object of this paper is to compute the Hochschild homology, the
cyclic homology and the Chern characters of the Auslander algebras $\gln$ of
the Taft algebras $\lan$ and of their Auslander algebras $\gln,$ in order to study a possible influence of  the Hopf algebra structure of $\lan$ on them.

Note that Auslander algebras are useful when considering Artin algebras of finite representation type, since there is a bijection between the Morita equivalence classes of such algebras  and the Morita equivalence classes of Auslander algebras (\emph{cf.}~\cite{ARS}).

The Hopf algebra structure on an algebra $\Lambda$ conveys an
additional structure on the Grothendieck groups $\K_0(\Lambda)$ and
$\overline{\K}_0(\Lambda) $ of isomorphism classes of projective
(\emph{resp.} all) indecomposable modules, since the tensor product
over the base ring $k$ of two $\Lambda-$modules is again a
$\Lambda-$module,\emph{via} the comultiplication of $\Lambda.$
Furthermore, there is a one-to-one correspondence between the
indecomposable modules over any algebra and the indecomposable
projective modules over its Auslander algebra; in the case of a Hopf
algebra, therefore, the Grothendieck group of projective modules of
$\gl$ is endowed with a multiplicative structure. However, this
correspondence does not preserve the underlying vector spaces, and
this multiplicative structure does not appear to be natural.

In this paper, we study the example of the Taft algebras; they are Hopf algebras which are neither commutative, nor cocommutative. They are interesting for various reasons; for instance, $\Lambda _p$ is an example of a non-semisimple Hopf algebra whose dimension is the square of a prime (\emph{cf.}~\cite{M1}). They are of finite representation type; furthermore, when $n$ is odd, $\lan$ is isomorphic to the half-quantum group \uq \hspace{1mm} ($q$ primitive $n^{th}-$root of unity), and is the only half-quantum group $u_q^+(\frak{g})$ at a root of unity which is not of wild representation type (\emph{cf.}~\cite{C1}). Then, for each $n,$ $\lan$ is not braided, but its Grothendieck group is a commutative ring nonetheless (\emph{cf.}~\cite{C2,G}). 

These examples show that the non-commutative, non-cocommutative Hopf algebra structure of $\lan$ does not yield a natural multiplicative structure on its cyclic homology. There is a product, however, obtained by transporting that of $\K_0(\Lambda)$ via the Chern characters, which are onto.

The paper is organised as follows: first, we recall the quiver of the
Auslander algebras $\gln$ and give a minimal projective resolution of
$\gln$ as a $\gln$-bimodule. Then, we compute the Hochschild and
cyclic homologies of these algebras, and finally we compute the Chern
characters of the algebras $\lan$ and  $\gln.$

Throughout this text, $k$ is an algebraically closed field.

\section{The quivers of the Auslander algebras $\gln$}

In this paragraph, we shall describe the objects of our study: the
Auslander algebras of the Taft algebras.

The Taft algebra $\lan$ is described by quiver and relations as
follows: the quiver is an oriented cycle with $n$ vertices and $n$
arrows, and the relations are all the paths of length greater than or
equal to $n$.

Its Auslander algebra has been described in \cite{GR} (\emph{see}
\cite{ARS} p232 for a general definition). Its quiver is $$\epsfysize=7cm\epsfbox{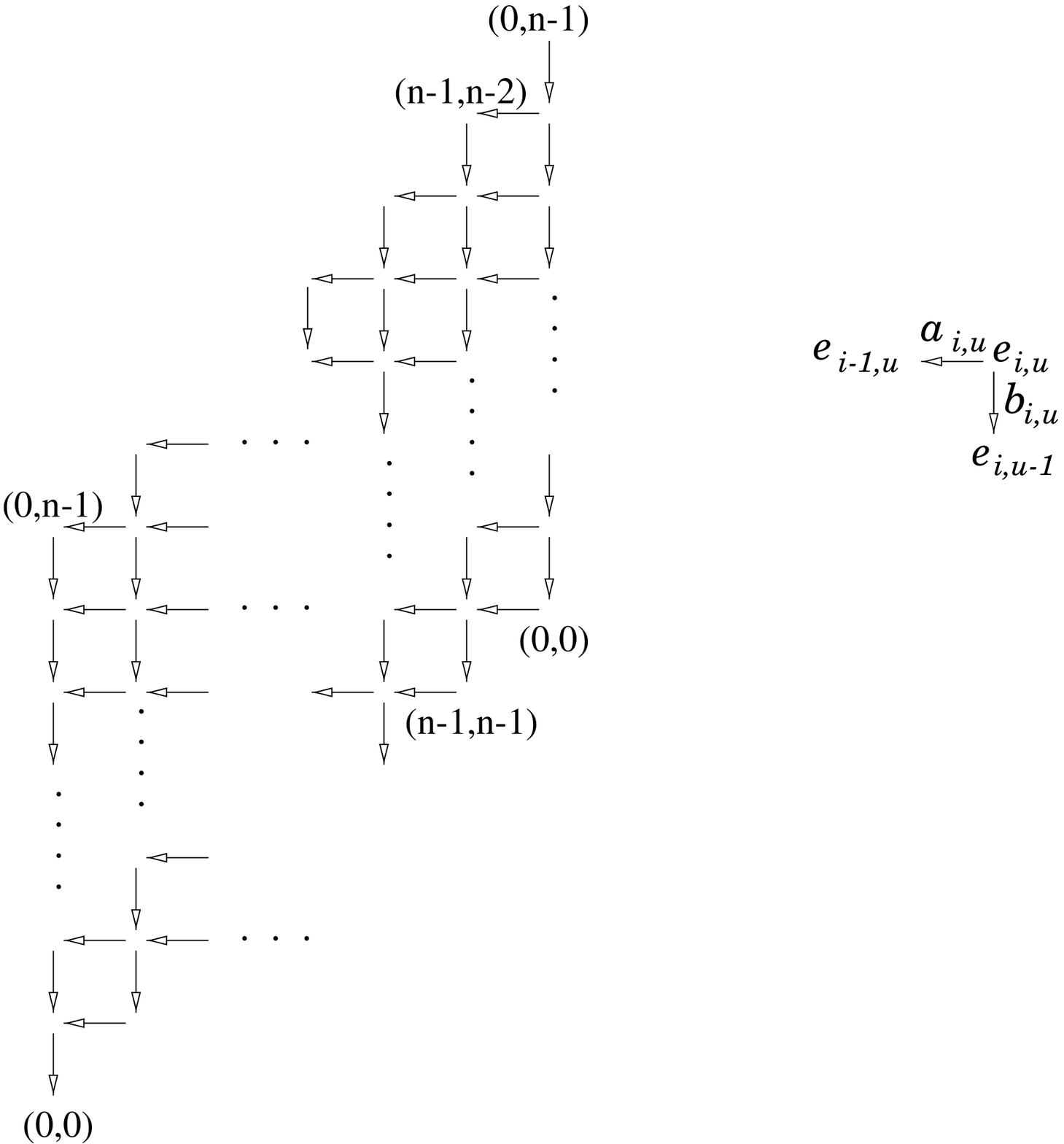}$$
where both vertical outer edges are identified (the quiver is on a cylinder). Let $Q_n$ denote this quiver,  let $\{e_{i,u}/(i,u) \in \Bbb{Z}/n\Bbb{Z} \times \Bbb{Z}/n\Bbb{Z}\}$ be the set of vertices of $Q_n,$ and let $\{a_{i,u}\, ; \,b_{i,u}/ \,(i,u) \in \Bbb{Z}/n\Bbb{Z} \times \Bbb{Z}/n\Bbb{Z}\}$ be the set of edges of $Q_n,$ as in the figure above. 

The mesh relations on this quiver are: $a_{i,i-2}\, b_{i,i-1}=0$ for all $i \in \Bbb{Z}/n\Bbb{Z}$ (the composition of two edges of any `triangle' under the top diagonal is zero), and $a_{i,u-1}b_{i,u}+b_{i-1,u}a_{i,u}=0$ for all $i$ and $u$ in $\Bbb{Z}/n\Bbb{Z}$ (the squares are anticommutative). 

The algebra $\gln$ is the quotient of the path algebra $kQ_n$ by the
ideal generated by these relations. We shall assume $n\geq 2$.


\section{Hochschild and cyclic homologies of $\gln$}

We are going to use the following theorem due to Happel to compute a
minimal projective resolution of $\gln$ as a $\gln$-bimodule (the
situation in \cite{H}  is more general):

\begin{theo}(\cite{H} 1.5)\label{happel} If $$\ldots \rightarrow R_p
  \rightarrow R_{p-1} \rightarrow \ldots \rightarrow R_1 \rightarrow
  R_0 \rightarrow \gln \rightarrow 0$$ is a minimal projective
  resolution of $\gln$ as a $\gln$-bimodule, then $$R_p=
  \bigoplus_{\tiny \begin{array}{c} (i,u) \\(j,v) \end{array}} (\gln
  e_{j,v} \ot e_{i,u} \gln )^{\mathrm{dim} _k Ext_{\gln}^p (S_{i,u} \,
    ; S_{j,v})}.$$ \end{theo}

In our case, we have:

\begin{cor} The complex
$$ \begin{array}{rcl} \ldots 0 \rightarrow 0 \rightarrow \bigoplus_{\tiny
  \begin{cases}(i,u) \\ i \neq u \end{cases}} \gln e_{i-1,u-1} \ot
e_{i,u} \gln &\rightarrow& \bigoplus_{\tiny (i,u)} [ (\gln e_{i-1,u} \ot
e_{i,u} \gln) \oplus  (\gln e_{i,u-1} \ot e_{i,u} \gln)]\\&\rightarrow&
\bigoplus_{\tiny (i,u)} \gln e_{i,u} \ot e_{i,u} \gln  \rightarrow
\gln \rightarrow 0\end{array}$$ is a minimal projective resolution of $\gln$ as a
  $\gln$-bimodule. 
\end{cor}

\begin{proof}
We need to compute the $\Ext$ groups between the simple
$\gln$-modules. First, let us compute the projective resolutions of
the simple modules:

\begin{lemma} Let $P_{i,u}$ denote the indecomposable projective
  $\gln-$module at the vertex $e_{i,u},$ and let $S_{i,u} = \topm
  (P_{i,u})$ be the corresponding simple module. The minimal projective resolutions of the simple modules  are: $$\begin{array}{ccccccccr} 
&&0 \ \longrightarrow & P_{i-1,i} &  \longrightarrow & P_{i,i} & \longrightarrow & S_{i,i} &  \longrightarrow  0 \\
0  \longrightarrow & P_{i-1,i-2} & \longrightarrow & P_{i,i-2} & \longrightarrow & P_{i,i-1} & \longrightarrow & S_{i,i-1} & \longrightarrow 0 \\
0  \longrightarrow & P_{i-1,i-j-1} & \longrightarrow & P_{i-1,i-j} \oplus P_{i,i-j-1} & \longrightarrow & P_{i,i-j} & \longrightarrow & S_{i,i-j} & \longrightarrow 0
\end{array}$$ for $2 \leq j \leq n-1.$ \end{lemma}

\begin{prooflemma} We consider only the $S_{n-1,u},$ because the other
  cases may be obtained by translating the quiver along the cylinder
  on which it lies. It is then straightforward to compute their
  minimal projective resolutions. \end{prooflemma}

We can now compute the $\Ext$ groups:

\begin{lemma}  Let $S$ be a simple  $\gln$-module. Then:
\begin{eqnarray*}
\Ext ^0_{\gln} (S_{i,u}\, ; S)& = &\begin{cases} k & \mbox{ if $S=S_{i,u}$,} \\ 0 & \mbox{ if $S \neq S_{i,u}$,}\end{cases} \\
\Ext ^1_{\gln} (S_{i,i}\, ; S)& = &\begin{cases} k & \mbox{ if $S=S_{i-1,i}$,} \\ 0 & \mbox{ if $S \neq S_{i-1,i}$,}\end{cases}\\
\Ext ^1_{\gln} (S_{i,i-1}\, ; S)& =&\begin{cases} k & \mbox{ if $S=S_{i,i-2}$,} \\ 0 & \mbox{ if $S \neq S_{i,i-2}$,}\end{cases} \\
\Ext ^1_{\gln} (S_{i,i-j}\, ; S)& =&\begin{cases} k & \mbox{ if $S=S_{i,i-j-1}$ or  $S=S_{i-1,i-j}$,} \\ 0 & \mbox{ if $S \neq S_{i-1,i-j-1}$ and  $S \neq S_{i-1,i-j}$,}\end{cases} \mbox{  if 1 $\leq j \leq n-1$} \\
\Ext ^2_{\gln} (S_{i,i}\, ; S)& =& 0\\
\Ext ^2_{\gln} (S_{i,i-j}\, ; S)& =&\begin{cases} k & \mbox{ if $S=S_{i-1,i-j-1}$,} \\ 0 & \mbox{ if $S \neq S_{i-1,i-j-1}$,}\end{cases} \mbox{  if 1 $\leq j \leq n-1$} \\
\Ext ^p_{\gln} (S_{i,u}\, ; S)& = & 0 \mbox{ if $ p \geq 3.$}
\end{eqnarray*} \end{lemma}

Applying Happel's Theorem (\ref{happel}) we get the minimal projective resolution for
$\gln.$ \end{proof}

 Now applying the functor $\gln \ot _{\tiny \gln - \gln} ? \ $ to this resolution, we obtain a complex: $$\ldots 0 \longrightarrow \ldots \longrightarrow 0 \longrightarrow k Q_0 \longrightarrow 0.$$ Therefore:

\begin{prop} The Hochschild homology of $\gln$ is: $$\begin{cases} \Hh_0 (\gln) = kQ _0 \cong k^{n^2} \\
\Hh_p (\gln) =0 & \forall p >0, \end{cases}$$ and hence the cyclic homology of $\gln$ is: $$\begin{cases} \Hc_{2p} (\gln) = kQ _0 \cong k^{n^2} \\
\Hc_{2p+1} (\gln) =0  \end{cases}  \forall p \geq 0.$$
\end{prop}

\begin{rmk} There doesn't seem to be any connection between these
  results and those for $\lan$. Indeed, the Hochschild and cyclic
  homologies for the Taft algebras are given as follows (\emph{see}
  \cite{S} for the Hochschild homology and 
  \cite{T,T1} for the cyclic homology):
$$\begin{cases}\Hh_{0}  (\Lambda _n ) = k^n \\
\Hh_{p} (\Lambda _n ) =k^{n-1}& \forall p >0.
\end{cases}$$ and
$$\begin{cases} \Hc_{2p} (\Lambda _n ) = k^n, \\
\Hc_{2p+1} (\Lambda _n ) = k^{n-1} & \forall p \geq 0.
\end{cases}$$
\end{rmk}


\section{Chern characters of $\lan$ and $\gln$}

Let $\K_0(\lan)$ (\emph{resp.} $ \K_0(\gln)$) be the Grothendieck group of projective $\lan $-modules (\emph{resp.} $\gln $-modules). We are interested in the Chern characters $ch_{0,p} : \K_0(\lan) \rightarrow  \Hc_{2p} (\lan )$ (\emph{resp.} $ \K_0(\gln) \rightarrow \Hc_{2p} (\gln)$). We shall write $[P_j]$ (\emph{resp.} $[P_{i,u}]$) for the isomorphism class of the projective module at the vertex $e_j$ (\emph{resp.} $e_{i,u}$).

Set $\sigma ^p = (y_p,z_p,\ldots, y_1,z_1,y_0) \in \Bbb{N}^{2p+1}$
with $y_p=(-1)^p (2p)!/p!$ and $z_p=(-1)^{p-1} (2p)!/2(p!).$ There is
a system of generators of $\Hc_{2p} (\lan)$ (\emph{resp.} $\Hc_{2p} (\gln)$)
given by the following set: $$\{\sigma _i^p := \sigma ^p (e_i, \ldots
, e_i)\in (\Tot \, CC(\lan))_{2p}\; /\;i=0, \ldots, n-1 \}$$ (\emph{resp.} by
$\{\sigma _{i,u}^p := \sigma ^p (e_{i,u} , \ldots , e_{i,u})\in (\Tot \,  CC(\gln))_{2p}\; /\;i, u \in \{0,1, \ldots, n-1\}  \}).$

Consider the elements $$\begin{array}{rclcrcl}\epsilon_j : \lan &
  \longrightarrow & \lan &\mbox{and} & \epsilon_{i,u} : \gln & \longrightarrow & \gln \\ \lambda & \mapsto & \lambda e_j &&\lambda & \mapsto & \lambda e_{i,u}\end{array}$$ in $\mathcal{M}_1(\lan)$ and  $\mathcal{M}_1(\gln);$ their ranges are the corresponding projective modules. Then by definition of the Chern characters (\emph{see} \cite{L} 8.3.4), we have: $$\begin{array}{l} 
ch_{0,p} ([P_j]) = ch_{0,p} ([\epsilon_j]) := \mathrm{tr} (c(\epsilon_j))= \sigma _j^p  \ \mbox{in $\Hc_{2p}(\lan)$}\\
ch_{0,p} ([P_{i,u}])= ch_{0,p} ([\epsilon_{i,u}]) = \sigma _{i,u}^p  \ \mbox{in $\Hc_{2p}(\gln)$}.
\end{array}$$ using the isomorphisms $\mathcal{M}_m(\Lambda) \cong \mathcal{M}_m(k) \ot \Lambda.$ Here, $$c(\epsilon_j)=(y_p \epsilon_j ^{\ot 2p+1} , z_p \epsilon_j ^{\ot 2p}, \ldots , z_1 \epsilon_j ^{\ot 2}, y_0 \epsilon_j) \in \mathcal{M}(\gln)^{\ot 2p+1} \oplus \ldots \oplus \mathcal{M}(\gln).$$
 
\begin{rmk} There is a decomposition formula for the tensor product of indecomposable modules on $\lan$ (\emph{see}~\cite{C2,G}). From this formula, we get inductively: $$ch_{0,p} ([L_1] \ot \ldots \ot [L_r]) = \frac{1}{n^2} \prod_{i=1}^r \left( \mathrm{dim}\, L_i \right) \; (\sigma _0^p ,\ldots , \sigma _{n-1} ^p), \ \mbox{for $r\geq 2,$}$$ where the $L_i$ are arbitrary projective $\lan$-modules. Unfortunately, this product in the cyclic homology doesn't seem natural.
\end{rmk}

\begin{rmk} Let $\overline{\K}_0 (\lan)$ be the Grothendieck group of all $\lan$-modules (not just the projective ones). Then $\overline{\K}_0 (\lan) \cong \K_0(\gln).$ Hence, if $N_{i,u} $ is the indecomposable $\lan$-module which starts at the vertex $i$ and ends at the vertex $u,$ it corresponds to the projective $\gln$-module $P_{i,u},$ and we get a map: 
\begin{eqnarray*} 
\overline{\K}_0 (\lan)& \longrightarrow & \Hc_{2p} ( \gln) \\ 
N_{i,u} & \mapsto & \sigma_{i,u}^p. 
\end{eqnarray*}
\end{rmk}

\begin{rmk} Although $\gln$ is not a Hopf algebra, its Grothendieck
  group $\K_0(\gln)$ does have a ring structure, which does not appear
  to be natural: for every $[P]$ in $\K_0(\gln),$ there exists a $[B]$ in $\K_0(\lan)$ such that $[P]=[\Hom_{\lan}(M,B)],$ where $M$ is the sum of all isomorphism classes of indecomposable $\lan$-modules. If $[Q]=[\Hom_{\lan}(M,C)]$ is another element in  $\K_0(\gln),$ we can set $$[P].[Q]= [\Hom_{\lan}(M,B\ot _k C)]$$ (the vector space $B\ot _k C$ is a $\lan$-module since $\lan$ is a Hopf algebra).  In fact, using the decomposition in \cite{C2,G}, the product can be written: $$[P_{i,u}][P_{j,v}]=  \begin{cases} \sum_{l=0}^{v-j} [P_{i+j+l, u+v -l}] & \mbox{ if $u+v-(i+j) \leq n-1$}\\ \sum_{l=0}^{e} [P_{i+j+l, u+v +l-1}] + \sum_{m=e+1}^{v-j} [P_{i+j+m, u+v -m}] & \mbox{ if $e:=u+v-(i+j) -( n-1) \geq 0$}
\end{cases}$$  where $u$ and $v$ represent
elements in $\mathbb{Z}/n\mathbb{Z}$ such that $u-i$ and $v-j$ are in $\{ 0, 1
, \ldots , n-1 \}.$
\end{rmk}


\flushleft{{\hrulefill\hspace*{10cm}}\\
\footnotesize{ 
{\sc Rachel Taillefer\\Department of Mathematics and Computer
    Science,\\ University of Leicester,\\ Leicester LE1 7RH, \\United
    Kingdom.\\ E-mail:} R.Taillefer@mcs.le.ac.uk.}}

\end{document}